\documentclass[11pt]{amsart}
\usepackage{epsfig}
\usepackage{amsthm}
\usepackage{amssymb}
\usepackage{amsmath}
\usepackage{epic}
\usepackage{eepic}

\newcommand     {\comment}[1]   {}
\newcommand{\mute}[2] {}
\newcommand     {\printname}[1] {}
\newcommand{\labell}[1] {\label{#1}\printname{#1}}

\newtheorem {Theorem}   {Theorem}
\newtheorem {Claim}[equation]    {Claim}
\newtheorem {Proposition}[equation]{Proposition}
\theoremstyle{definition}

\theoremstyle{remark}
\newtheorem{Remark}[equation]{Remark}

\newtheorem {Corollary}[equation]{Corollary}

\def    \valpha      {\alpha}
\def    \vV {{V}}

\def    \bfq    {{\mathbf q}}

\def    \bfC    {{\mathbf C}}

\def    \bfO    {{\mathbf O}}

\def    \bfN    {{{\mathbf N}_q}}

\def    \inv    {^{-1}}

\def    \ssminus        {{\smallsetminus}}

\def    \l<        {\left< }
\def    \r>        {\right> }
\def    \del    {\partial}
\newcommand {\deldel}[1] {\frac{\partial}{\partial #1}}

\def    \nonneg {{\geq 0}}

\def    \span   {{\operatorname{span}}}

\def    \Z  {{\mathbb Z}}
\def    \R  {{\mathbb R}}
\def    \C  {{\mathbb C}}

\newcommand{\e}{\mathrm{e}}       % exponential map
\newcommand{\dsum}{\displaystyle\sum}  % Sum in displaystyle
\newcommand{\dint}{\displaystyle\int}  % Integral in dysplaystyle

\newcommand{\lclass}{\mathbf{L}}  % L-class
\newcommand{\todd}{\mathbf{Td}}   % Todd-class

\begin{document}

\bibliographystyle{plain}

\title[The Weighted Euler-Maclaurin Formula for a simple integral polytope]
{The Weighted Euler-Maclaurin Formula for a simple integral
polytope}

\author[J.\ Agapito]{Jos\'e Agapito}
\address{Department of Mathematics, University of California,
Santa Cruz, CA 95064, USA} \email{jarpepe@math.UCSC.EDU}

\author[J.\ Weitsman]{Jonathan Weitsman}\thanks{This work was partially supported
by National Science Foundation Grant DMS 99/71914.}
\address{Department of Mathematics, University of California,
Santa Cruz, CA 95064, USA}  \email{weitsman@math.UCSC.EDU}
\thanks{2000 \emph{Mathematics Subject Classification.}
Primary 65B15, 52B20.}

\begin{abstract}
We give an Euler-Maclaurin formula with remainder for the weighted
sum of the values of a smooth function on the integral points in a
simple integral polytope. Our work generalizes the formula
obtained in~\cite{ArX}.
\end{abstract}
\maketitle
%\tableofcontents

% -------------------------------------------------------------------------
\section{Introduction}
% -------------------------------------------------------------------------
\labell{sec:intro}

The Euler-Maclaurin formula computes the sum of the values of a
function $f$ over the integer points in an interval in terms of
the integral of $f$ over variations of that interval. Khovanskii
and Pukhlikov \cite{KP1,KP2} and Kantor and Khovanskii \cite{KK}
generalized the classical Euler-Maclaurin formula to higher
dimensional convex polytopes $\Delta$ which are integral and
regular. This formula was generalized to \emph{simple} integral
polytopes by Cappell and Shaneson \cite{CS:bulletin,CS:EM}, and
subsequently by Guillemin \cite{G} and by Brion-Vergne \cite{BV}.
These generalizations involve corrections to Khovanskii's formula
when the simple polytope is not regular. These formulas are exact
formulas, valid for sums of exponential or polynomial functions.
With the use of the $\lclass$ class associated to the signature
operator, as in \cite{CS:bulletin,CS:EM}, Karshon, Sternberg and
Weitsman gave an Euler-Maclaurin formula with remainder for a
weighted sum of the values of an arbitrary smooth function on the
lattice points in a simple integral polytope \cite{PNAS,ArX}.

The purpose of this paper is to give a generalization of the
Euler-Maclaurin formula with remainder of \cite{ArX}, to allow for
more general weightings, including the ordinary, unweighted sum.
To do this, we use the Hirzebruch formal power series
${\chi}_q(S)$ \cite{H} and the weighted polar decomposition of
\cite{A}.

% -------------------------------------------------------------------------
\section{Weighted Euler-Maclaurin in one dimension}
% -------------------------------------------------------------------------
\labell{sec:wem-1d}

\subsection*{Weighted sums in one dimension}

Let $q$ be any complex number and $f(x)$ be any function on the
real line. For any integers $a<b$, define
\begin{equation} \labell{wted sum interval}
 {\sum_{[a,b]}}^q f := q f(a) + f(a+1) + \ldots + f(b-1) +
 q f(b) .
\end{equation}
Similarly, for a ray $[a,\infty)$ with $f$ having compact support
\begin{equation}\labell{wted sum ray}
{\sum_{[a,\infty)}}^q f := qf(a)+f(a+1)+\ldots
\end{equation}
so
\begin{equation}\label{1dweighedLawrence}
{\sum_{[a,b]}}^q f = {\sum_{[a,\infty)}}^q f -
{\sum_{[b,\infty)}}^{1-q} f .
\end{equation}

If $q=1$ this is the ordinary sum; if $q=1/2$ this is the weighted
sum of \cite{ArX}.

\subsection*{Euler-Maclaurin formulas}

The classical Todd function is defined by
\begin{equation*}
\todd(S)=\dfrac{S}{1-\e^{-S}} = 1 - b_1 S + \dsum_{n=1}^\infty
\frac{b_{2n}}{(2n)!} S^{2n} ,
\end{equation*}
where $b_n$ is the $n$-th Bernoulli number (we are using the
conventions in \cite{Bo}). Similarly, the $\bf L$ function is
given by
\begin{equation*}
\lclass(S)=\dfrac{S/2}{\tanh(S/2)} = 1 + \dsum_{n=1}^\infty
\frac{b_{2n}}{(2n)!}S^{2n} .
\end{equation*}
Both $\todd(S)$ and $\lclass(S)$ are convergent power series for
$|S| < 2 \pi.$\par

The Hirzebruch function is defined by
\begin{equation}\label{XS}
{\bf\chi}_q(S) = q\todd(S) + (1-q)\todd(-S) = 1 + (q-\frac{1}{2})S
+ \dsum_{n=1}^\infty \frac{b_{2n}}{(2n)!}S^{2n} .
\end{equation}

It is related to $\todd$ and to $\lclass$ by
\begin{equation}\labell{XS-Td-L}
{\bf\chi}_q(S) = (q-1)S + \todd(S) = \left(q-\frac{1}{2}\right) S
+ \lclass(S) .
\end{equation}

Thus if $q=1$ we have ${\bf\chi}_1(S)=\todd(S)$ and if $q=1/2$, we get
${\bf\chi}_{1/2}(S)=\lclass(S)$.
\bigskip

Let $f$ be a compactly supported function on the real line of
class $C^m$ where $m>1$. The standard Euler-Maclaurin formula with
remainder for a ray can be written (see for instance \cite{Bo})
\begin{equation}\labell{st-emr1dray}
 {\sum_{[a,\infty)}} f
 = \left. \todd^{2k}(\deldel{h}) \int_{a-h}^{\infty} f(x) dx
 \right|_{h=0} + R_{m}^a(f),
\end{equation}
where $k = \lfloor m/2 \rfloor$, where $\todd^{2k}(S)$ denotes the
truncation of the power series $\todd(S)$ at the $2k$-th term, and
where the remainder $R_{m}^a(f)$ is given by
$$R_m^a(f) = (-1)^{m-1} \int_a^\infty P_m(x) f^{(m)}(x) dx\quad\mbox{with}
\quad P_m(x) = \frac{B_m(\{ x \})}{m!} \,\, .$$
Here, $B_m(x)$ is
the $m$th Bernoulli polynomial (see \cite{Bo}) and $\{ x \} = x -
\lfloor x \rfloor$ is the fractional part of $x$. Moreover, the
function $P_m(x)$ is given by
\begin{equation} \labell{P2k}
 P_{2k}(x) = (-1)^{k-1} \sum_{n=1}^\infty \frac{2\cos2\pi nx}{(2\pi n)^{2k}}
\end{equation}
if $m=2k$, and by
\begin{equation} \labell{P2kp1}
 P_{2k+1}(x) =
   (-1)^{k-1} \sum_{n=1}^\infty \frac{2\sin2\pi n x}{(2\pi n)^{2k+1}}
\end{equation}
if $m=2k+1$.\par

Similarly, denoting the truncation of the power series
${\bf\chi}_q$ at the $2k$-th term by ${\bf\chi}^{2k}_q$, we have a
weighted Euler-Maclaurin formula with remainder for a ray,
\begin{equation}\label{emr1dray+}
 {\sum_{[a,\infty)}}^q f
 = \left. \chi_q^{2k}(\deldel{h}) \int_{a-h}^{\infty} f(x) dx
 \right|_{h=0} + R_{m}^a(f),
\end{equation}

\noindent where $f$ is a compactly supported function.

Combining \eqref{emr1dray+} and \eqref{1dweighedLawrence} we get

\begin{equation}
{\sum_{[a,b]}}^q f = \left. \chi_q^{2k}(\deldel{h_1})
\int_{a-h_1}^{\infty} f(x) dx \right|_{h_1=0} - \left.
\chi_{1-q}^{2k}(\deldel{h_2}) \int_{b-h_2}^{\infty} f(x) dx
\right|_{h_2=0}\\
+ R_m^{[a,b]}(f).
\end{equation}

\noindent where

$$R_m^{[a,b]}(f) = R_m^b(f) - R_m^a(f)
= (-1)^{m-1} \int_a^b P_m(x) f^{(m)}(x) dx.$$

Given the fact that $\chi_q^{2k}(S)$ is a polynomial whose constant
term is $1$, we can write this as

\begin{equation}\label{emrint}
{\sum_{[a,b]}}^q f =
\chi_q^{2k}(\deldel{h_1})\chi_{1-q}^{2k}(\deldel{h_2})
\left.\left(\int_{a-h_1}^{\infty} f(x) dx -\int_{b-h_2}^{\infty}
f(x) dx\right) \right|_{h_1=h_2=0}\\
+ R_m^{[a,b]}(f).
\end{equation}

To make further progress, we note the following symmetry property.
\begin{equation}\label{symmetry}
{\bf\chi}^{2k}_q(S)={\bf\chi}^{2k}_{1-q}(-S).
\end{equation}

(To see this, observe that ${\bf\chi}^{2k}_q(S)$ is a polynomial
with constant coefficients whose constant term is 1, whose linear
term is $(q-1/2)S$ (see \eqref{XS-Td-L}), and whose other terms
are all of even degree independent of $q$.)

Thus we obtain the following result:

\begin{Theorem}[Euler-Maclaurin with remainder for intervals]
\label{EMintervalgrm} Let $f(x)$ be a function with $m > 1$
continuous derivatives and let $k = \lfloor m/2 \rfloor$. Then
\begin{equation}\label{emr1d}
  {\sum_{[a,b]} }^q f = \left. {\bf\chi}_q^{2k}(\deldel{h_1})
  {\bf\chi}_q^{2k}(\deldel{h_2})
  \int_{a-h_1}^{b+h_2} f(x) dx
  \right|_{h_1 = h_2 = 0} + R_m^{[a,b]}(f) \, .
\end{equation}
\end{Theorem}

Note that our argument applies to functions $f$ of compact support.
However, for a general function $f$ of type $C^m$, the theorem remains
true:  We need only multiply $f$ by a smooth function of
compact support which is identically one in a neighborhood of $[a,b].$

If $f$ is a polynomial, Theorem \ref{EMintervalgrm}
becomes exact when $m$ is greater than the degree of $f$:

\begin{Corollary} Let $f$ be a polynomial. Then
\begin{equation} \labell{EM exact interval}
   {\sum_{[a,b]}}^q f =
   \left. {\bf\chi}_q(\deldel{h_1}) {\bf\chi}_q(\deldel{h_2})
   \int_{a-h_1}^{b+h_2} f(x) dx \right|_{h_1 = h_2 = 0}.
\end{equation}
\end{Corollary}

\subsection*{Twisted weighted Euler-Maclaurin for a ray}

Consider the ``twisted weighted sum"
\begin{equation}\label{twistedsum-1dray}
{\sum_{n\ge 0}}^q \lambda^n f(n) = qf(0) + \sum_{n=1}^\infty
\lambda^n f(n),
\end{equation}
where $\lambda \neq 1$ is a root of unity, say, of order $N$.\par

Let $Q_{m,\lambda}(x)$ be distributions successively defined by
$$ Q_{0,\lambda}(x) = - \sum_{n \in \Z} \lambda^n \delta(x-n)$$
and
$$ \dfrac{d}{dx}Q_{m,\lambda}(x) =
Q_{m-1,\lambda}(x)\quad\mbox{and}\quad\dint_0^N Q_{m,\lambda}(x)
dx = 0 .$$

These distributions appear in \cite{ArX}.\par

Now, define the polynomial
$${\bf N}_q^{k,\lambda}(S) =
\left( q + \frac{\lambda}{1-\lambda} \right) S
 + Q_{2,\lambda}(0) S^2
 + Q_{3,\lambda}(0) S^3
 + \cdots
 + Q_{k,\lambda}(0) S^k,$$
for a root of unity $\lambda \neq 1$. When $q=1/2$, we get the
polynomial ${\bf M}^{k,\lambda}(S)$ defined in \cite{ArX}. Since
${\bf N}_q^{k,\lambda}(S)$ and ${\bf M}^{k,\lambda}(S)$ differ by
$(q-1/2)S$, adding this term gives the following generalization of
Proposition 23 in \cite{ArX}.

\begin{Proposition}\labell{twistedemray}
Let $k>1$ and let $f \in C_c^{k}(\R)$. Then
\begin{equation} \labell{twistedemrayeq}
\left.
  {\sum_{n\ge 0}}^q \lambda^n f(n) =
  {\bf N}_q^{k,\lambda}(\deldel{h}) \int_{-h}^\infty f(x) dx \right|_{h=0}
  + (-1)^{k-1}\int_0^\infty Q_{k,\lambda}(x) f^{(k)}(x) dx .
\end{equation}
\end{Proposition}

As in \cite{ArX}, we have the following symmetry property
\begin{equation}\labell{Dandlambdainverse}
{\bf N}_{1-q}^{m,1/\lambda}(S)= {\bf N}_{q}^{m,\lambda}(-S).
\end{equation}

\begin{Remark} \label{M vs L}
For $\lambda = 1$, we define
$$ {\bf N}_q^{k,1}(S) = \chi_q^{2 \lfloor k/2 \rfloor}(S) \quad \text{and} \quad
   Q_{k,1} = P_k,$$
then \eqref{twistedemrayeq} boils down to \eqref{emr1dray+}. So,
\eqref{twistedemrayeq} also holds for $\lambda = 1$. Notice that
if $\lambda \neq 1$ then ${\bf N}_q^{k,\lambda}(S)$ is a multiple
of $S$, and that if $\lambda = 1$ then ${\bf N}_q^{k,\lambda}(S) =
1 + $ a multiple of $S$. Property \eqref{Dandlambdainverse}
continues to hold for $\lambda = 1$ because of the symmetry
property \eqref{symmetry}. Finally, if $\lambda=1$, we get the
truncation $\todd^{k}(S)$ of $\todd(S)$ for $q=1$, and the
corresponding truncation at $k$ of $\todd(-S)$ for $q=0$. These two expressions
differ by $S$.
\end{Remark}

% -----------------------------------------------------------------------------
\section{Simple polytopes and finite groups associated to them}
% -----------------------------------------------------------------------------

In this section we recall various combinatorial and group-theoretic
facts about simple polytopes which will be needed in our proof of
the weighted Euler-Maclaurin formula with remainder.  Most of this
material is taken from \cite{A} and \cite{ArX} and is included for completeness.
\labell{sec:prelim}

Let $H_i = \{ x \mid \langle u_i , x \rangle + \mu_i \geq 0 \,,\, \mu_i\in\R \}$ be
a half space in $\R^n$, where $u_i\in{\R^n}^*$ for $i=1,\ldots,d$. A compact convex
polytope $\Delta$ in $\R^n$ is a compact set which can be written as the
intersection of finitely many half-spaces
\begin{equation} \labell{Delta}
\Delta = H_1 \cap \ldots \cap H_d,
\end{equation}
with the smallest possible $d$, so that the $H_i$'s are uniquely
determined up to permutation. We order them arbitrarily. The
\emph{facets} (codimension one faces) of $\Delta$ are
$$ \sigma_i = \Delta \cap \del H_i \quad , \quad i = 1, \ldots, d .$$

The vertices of $\Delta$ are all possible codimension n faces
obtained by intersections of facets.\par

A polytope is called \emph{integral} if its vertices are in the
lattice $\Z^n$; it is called \emph{simple} if exactly $n$ edges
emanate from each vertex; it is called \emph{regular} if,
additionally, the edges emanating from each vertex lie along lines
which are generated by a $\Z$-basis of the lattice $\Z^n$.\par

For each vertex $v$ of $\Delta$, let $I_v \subset \{ 1, \ldots, d
\}$ encode the set of facets that contain $v$, so that
$$ i \in I_v \quad \text{if and only if} \quad v \in \sigma_i .$$

The vector $u_i \in {\R^n}^*$ can be thought of as the inward
normal to the $i$th facet of $\Delta$; a-priori it is determined
up to multiplication by a positive number. If the polytope
$\Delta$ is integral, then the $u_i$'s can be chosen to belong to
the dual lattice ${\Z^n}^*$, and we can fix our choice of the
$u_i$'s by imposing the normalization condition that the $u_i$'s
be primitive lattice elements, that is, that no $u_i$ can be
expressed as a multiple of a lattice element by an integer greater
than one.\par

Assume that $\Delta$ is \emph{simple}, so that each vertex is the
intersection of exactly $n$ facets. For each $i \in I_v$, there
exists a unique edge at $v$ which does not belong to the facet
$\sigma_i$; choose any vector $\alpha_{i,v}$ in the direction of
this edge. These vectors form a dual basis to the inward normal
vectors $u_i$'s and are uniquely determined when the polytope is
integral and the $u_i$'s are normalized in the sense explained
above.\par

A ``polarizing vector'' $\xi \in {\R^n}^*$ is a vector such that
$\langle \xi , \alpha_{i,v} \rangle$ is non-zero for all vertices
$v$ and all edges $i$ emanating from $v$. A choice of a polarizing
vector $\xi$ determine \emph{polarized edge vectors}
$\valpha^{\sharp}_{i,v}$ defined by
\begin{equation} \labell{flipalpha}
  \valpha^{\sharp}_{i,v} = \begin{cases}
 \valpha_{i,v} & \text{ if } \l< \xi , \valpha_{i,v} \r> < 0,\quad
  (\mbox{unflipped}) \\
 - \valpha_{i,v} & \text{ if } \l< \xi , \valpha_{i,v} \r> > 0,\quad
  (\mbox{flipped})
 \end{cases} .
\end{equation}

Let $q$ be any complex number. For each $i\in I_v$ define
\begin{equation}\label{flipq}
q_{i,v}^\sharp = \left\{%
\begin{array}{cl}
  q   & \mbox{if } \valpha^\sharp_{i,v} = \alpha_{i,v} \\
  1-q & \mbox{if } \valpha^\sharp_{i,v} = -\alpha_{i,v} \\
\end{array}%
\right. .
\end{equation}

The \emph{tangent cone} to $\Delta$ at $v$ is
\begin{equation}\label{tan cone}
 \bfC_v =
    \{ v + r (x-v) \ \mid \ r \geq 0 \, , \, x \in \Delta \}
 = v + \sum_{i \in I_v} \R_{\geq 0} \alpha_{i,v}.
\end{equation}

Similarly\, the \emph{polarized tangent cone} at $v$ is defined by
\begin{equation} \labell{pol tan cone}
 \bfC_v^\sharp = v + \sum_{i \in I_v} \R_{\geq 0} \alpha_{i,v}^\sharp.
\end{equation}

A simple integral orthant $\bfC$ in $\R^n$ is the intersection of $n$ half-planes in
general position,
$$\bfC = H_1\cap\ldots\cap H_n\quad\mbox{with}\quad H_i = \{x \mid
\langle u_i,x\rangle + \mu_i\ge 0 \,,\, \mu_i\in\R\}\quad\mbox{for}\quad
i=1,\ldots,n , $$ where the $u_i$'s are inward normals to the facets $\sigma_i =
\bfC \cap \del H_i$ of $\bfC$, which can be chosen to be primitive elements of the
dual lattice ${\Z^n}^*$, and whose vertex $v=\cap_{i=1}^n\{x \mid \langle
u_i,x\rangle + \mu_i = 0 \,,\, \mu_i\in\R\}$ is in $\Z^n$. This implies that
$\mu_i\in\Z$ for all $1\le i\le n$. If $\alpha_1 , \ldots, \alpha_n$ is the dual
basis to the $u_i$'s, that is,
$$ \l< u_j , \alpha_i \r> = \begin{cases}
1 & j=i \\ 0 & j \neq i ,
\end{cases}$$
\noindent then
\begin{equation}\label{simp-inte-orth}
\bfC = v + \sum_{j=1}^n \R_{\nonneg} \alpha_j.
\end{equation}

We associate a complex number $q_i$ to each facet $\sigma_i$ of $\bfC$, and define
the weighting function
\begin{equation}\label{weighting}
w (x)= \prod_{i\in I_F} q_i ,
\end{equation}
where $I_F$ denotes the set of facets in $\bfC$ meeting at the face $F$ (the
smallest dimensional face in $\bfC$ containing $x$.) If $x$ is in the interior of
$\bfC$, we set $w(x) = 1 $, and if $x\notin\bfC$, we set $w(x) = 0$.

For a simple integral orthant, we consider the weighted sum
\begin{equation}\label{weightedsum-cone}
{\sum_{\bfC \cap \Z^n}}^{\bf q} f = \sum_{x \in \bfC \cap
   \Z^n} w(x) f(x) ,
\end{equation}

\noindent where $\bf q$ denotes the $n$-tuple $(q_1,\ldots q_n)$ used in the
definition of $w(x)$.\par

Let $q$ be a complex number. For a polytope, we consider the weighted sum
\begin{equation}
{\sum_{\Delta \cap \Z^n}}^{q} f = \sum_{x \in \Delta \cap \Z^n}
   q^{c(x)} f(x) ,
\end{equation}\label{weightedsum-polytope}
\noindent where $c(x)$ is the codimension of the smallest dimensional face in
$\Delta$ containing $x$.

Given $q\in\C$, a cone $\bfC_v$ and a polarizing vector
$\xi\in{\R^n}^*$, we get the $n$-tuple
$\bfq_v^\sharp := (q_{1,v}^\sharp,\ldots, q_{n,v}^\sharp)$
where $q_{i,v}^\sharp$ is defined in \eqref{flipq}.
With this notation, the
weighted polar decomposition of \cite{A} is as follows.

\begin{Theorem}[\cite{A}] For any polarizing vector $\xi$, we
have
\begin{equation} \labell{wheee}
 {\sum_{\Delta \cap \Z^n}}^{q} f = \sum_v (-1)^{\# v}
 {\sum_{\bfC_v^\sharp \cap \Z^n}}^{\bfq^\sharp_v} f ,
\end{equation}
where we sum over the vertices $v$ of $\Delta$ and where $\# v$ is
the number of edge vectors at $v$ that are ``flipped'' by the
polarization process~(\ref{flipalpha}).
\end{Theorem}

In obtaining a weighted Euler-Maclaurin formula with remainder for
simple integral polytopes, we associate certain finite groups to
them. We recall some definitions and results from Section \S 5 in
\cite{ArX}.\par

Let us consider a simple integral orthant $\bfC$ with vertex $v\in\Z^n$, as given in
\eqref{simp-inte-orth}. To $\bfC$ we associate the finite group
\begin{equation} \labell{Gamma}
 \Gamma := {\Z^n}^* / \sum \Z u_i .
\end{equation}
This group is trivial exactly if $\bfC$ is regular.\par

Now, let $\Delta$ be a simple integral polytope in $\R^n$. For any
face $F$ of $\Delta$, let $I_F$ denote the set of facets of
$\Delta$ which meet at $F$. Because $\Delta$ is simple, the
vectors $u_i$, for $i \in I_F$, are linearly independent. Let $N_F
\subseteq {\R^n}^*$ be the subspace
$$ N_F = \span \{ u_i \mid i \in I_F \}.$$

To each face $F$ of $\Delta$ we associate a finite abelian group
$\Gamma_F$. Explicitly, the lattice
$$\vV_F = \sum_{i \in I_F} \Z u_i \, \subset \, N_F $$
is a sublattice of $N_F \cap {\Z^n}^*$ of finite index, and the
finite abelian group associated to the face $F$ is the quotient
\begin{equation} \labell{def GammaF}
\Gamma_F := (N_F \cap {\Z^n}^*) / \vV_F.
\end{equation}
If $F=v$ is a vertex, this is the same as the finite abelian group
associated to the tangent cone $\bfC_v$ as in \eqref{Gamma}.

Let $E$ and $F$ be two faces of $\Delta$ with $F \subseteq E$.
This inclusion implies that $I_E \subseteq I_F$, and hence
$\Gamma_E \subseteq \Gamma_F$.

We define a subset $\Gamma_F^\flat$ of $\Gamma_F$ by
\begin{equation} \labell{def:GammaFsharp}
\Gamma_F^\flat
  := \Gamma_F \ssminus
\bigcup_{\text{faces } E \text{ such that }  E \supsetneq F}
\Gamma_E.
\end{equation}
Then
\begin{equation} \labell{Gammas}
 \Gamma_v = \bigsqcup_{\{F:v \in F\}}  \Gamma_F^\flat.
\end{equation}

The map
\begin{equation} \labell{lambda j v}
 \lambda_{\gamma,j,v} :=
   e^{ 2 \pi i \left< \gamma ,\alpha_{j,v} \right> } ,
   \quad \text{ for } \gamma \in \Gamma_v
   \quad \text{ and } j \in I_v,
\end{equation}
is a well defined character and it is a root of unity.\par

\begin{Claim}[\cite{ArX}, Claim 61] \labell{claim1}
If $\gamma \in \Gamma_F$ and $j \in I_F$, then
$\lambda_{\gamma,j,v}$ is the same for all $v \in F$.
\end{Claim}

This allows us to define $\lambda_{\gamma,j,F}$ for $\gamma \in
\Gamma_F$ and $j \in I_F$ such that
$$ \lambda_{\gamma,j,F} = \lambda_{\gamma,j,v}
\quad \text{ for } \gamma \in \Gamma_F \text{ and } j \in I_F,
\text{ if } v \in F .$$

\begin{Claim}[\cite{ArX}, Claim 62]\labell{claim2}
If $\gamma \in \Gamma_F$ and $j \in I_v \ssminus I_F$ then
$\lambda_{\gamma,j,v}$ is equal to one.
\end{Claim}

This allows us to define $\lambda_{\gamma,j,F} = 1$ when $\gamma
\in \Gamma_F$ and $j \in \{ 1 , \ldots, d \} \ssminus I_F $. Then
\begin{equation} \labell{horse}
 \lambda_{\gamma,j,F} = \lambda_{\gamma,j,v}
\quad \text{ for } \gamma \in \Gamma_F \text{ and } 1 \leq j \leq
d , \text{ if } v \in F
\end{equation}
and
\begin{equation} \labell{lambda is one}
\lambda_{\gamma,j,F} = 1 \quad \text{ for } \gamma \in \Gamma_F
\text{ if } j \not\in I_F.
\end{equation}

\begin{Claim}[\cite{ArX}, Claim 65]\labell{claim3}
If $\gamma \in \Gamma_F^\flat$ and $j \in I_F$, then $
\lambda_{\gamma,j,F} \neq 1$.
\end{Claim}

% -----------------------------------------------------------------------------
\section{Weighted Euler-Maclaurin with remainder for simple integral
polytopes.}
% -----------------------------------------------------------------------------
\labell{sec:EM simple}

Let $\Gamma$ be the finite group \eqref{Gamma} associated to a simple integral
orthant $\bfC$ (as defined in \eqref{simp-inte-orth}) with vertex at $v\in\Z^n$. The
map $\gamma \mapsto \e^{ \left< \gamma, x \right> }$ defines a character on $\Gamma$
(\cite{ArX}, Lemma 52) whenever $x \in \sum \Z \alpha_j$, and this character is
trivial if and only if $x \in \Z^n$.\par

By a theorem of Frobenius,
$$ \frac{1}{| \Gamma |} \sum_{\gamma \in \Gamma}
\e^{ 2 \pi i \left< \gamma, x \right> } = \begin{cases}
 1 & \text{if } x \in \Z^n \\
 0 & \text{if } x \not \in \Z^n \\
\end{cases} $$
for all $x \in \sum \Z \alpha_j$. Then, for any vector
$\bfq=(q_1,\ldots,q_n)\in\C^n$ and any function $f(x)$ compactly
supported on $\R^n$,
\begin{eqnarray}
\nonumber
 {\sum_{\bfC \cap \Z^n}}^{\bfq} f
 & = & {\sum_x}^{\bfq} \left( \frac{1}{|\Gamma|}
\sum_{\gamma \in \Gamma} \e^{2\pi i \left< \gamma , x \right>} \right) f(x) \\
\labell{eqAA}
 & = & \frac{1}{|\Gamma|} \sum_{\gamma \in \Gamma} {\sum_x}^{\bfq}
       \e^{2\pi i \left< \gamma , x \right>} f(x)
\end{eqnarray}
where we sum over all
\begin{equation} \labell{x}
 x = v + m_1 \alpha_1 + \ldots + m_n \alpha_n ,
\end{equation}
with the $m_i$'s being non-negative integers.

The simple integral orthant $\bfC$ is the image of the standard
orthant $\bfO=\prod_{i=1}^n \R_{\ge 0}$ in $\R^n$ under the affine
map
$$ (t_1,\ldots,t_n) \mapsto v + \sum t_i \alpha_i.$$
This map sends the lattice $\Z^n$ onto the lattice $\sum \Z
\alpha_j$. Let us concentrate on one element $\gamma \in \Gamma$.
Because $v \in \Z^n$, from \eqref{x} we get
$$ \e^{2\pi i \left< \gamma , x \right>}
   = \prod_{j=1}^n \lambda_{j} ^{m_j}
\quad \text{where} \quad \lambda_{j} = \e^{2\pi i\left< \gamma ,
\alpha_j \right> } ,$$ so that the inner sum in (\ref{eqAA})
becomes
\begin{equation} \labell{eqBB}
{\sum_x}^{\bfq} \e^{2\pi i \left< \gamma , x \right>} f(x) =
{\sum_{m_1 \geq 0}}^{q_1} \lambda_{1}^{m_1} \cdots {\sum_{m_n \geq
0}}^{q_n} \lambda_{n}^{m_n}\,\, g(m_1,\ldots,m_n),
\end{equation}
where
$$ g(m_1,\ldots,m_n) =
   f(v + m_1 \alpha_1 + \ldots + m_n \alpha_n).  $$

Given $q\in\C$ and $k>1$, we had the twisted remainder formula
(see \eqref{twistedemrayeq})
$$ {\sum_{m \geq 0}}^q \lambda^m g(m)
 = \left. \bfN^{k,\lambda}(\deldel{h}) \int_{-h}^\infty g(t) dt \right|_{h=0}
 + (-1)^{k-1} \int_0^\infty Q_{k,\lambda}(t) g^{(k)}(t) dt$$
for all compactly supported functions $g(x)$ of type $C^k$, where
$k > 1$, where $\lambda$ is a root of unity, and where
$N_q^{k,\lambda}$ is a polynomial of degree $\leq k$.

Iterating this formula, the sum in \eqref{eqBB} can be written as
\begin{multline*}
  \begin{aligned}
  {\bf N}_{q_1}^{k,\lambda_{1}}(\deldel{h_1}) \int_{-h_1}^\infty \cdots
     {\bf N}_{q_n}^{k,\lambda_{n}}(\deldel{h_n}) \int_{-h_n}^\infty
     g(t_1,\ldots,t_n) dt_1 \cdots dt_n \\
    + R^{st}_{\bfq,k}(\lambda_{1},\dots,\lambda_{n};g)
  \end{aligned}                                 \\
     = \prod_{i=1}^n {\bf N}_{q_i}^{k,\lambda_i}(\deldel{h_i})
       \int\limits_{{\bf O}(h)} g(t_1,\ldots,t_n) dt_1 \cdots dt_n
      + R^{st}_{\bfq,k}(\lambda_{1},\ldots,\lambda_{n};g),
\end{multline*}
where
$$ {\bf O}(h_1,\ldots,h_n) = \{ (t_1,\ldots,t_n) \mid t_i \geq -h_i
\text{\ for all $i$} \}$$
denotes the ``dilated'' standard orthant, and where the remainder
is given by
\begin{multline} \labell{remainder}
R^{st}_{{\bfq},k}(\lambda_1,\dots,\lambda_n ; g) :=
\sum_{I \subsetneq \{ 1, \ldots, n \} } (-1)^{(k-1)(n-|I|)} \\
\prod_{i \in I} {\bf N}_{q_i}^{k,\lambda_i}(\deldel{h_i})
\int\limits_{{\bf O}(h)} \left. \prod_{i \notin I} Q_{k,\lambda_i}
(t_j) \prod_{i \notin I} \frac{\del^{k}}{\del t_i ^{k}}
 g(t_1,\ldots,t_n) dt_1 \cdots dt_n
\right|_{h=0}
\end{multline}
with
$$ g(t_1,\ldots,t_n) =
f(v+ t_1 \alpha_1 + \ldots + t_n \alpha_n).$$

Performing the change of variable given by the transformation
$$ L \colon (t_1 , \ldots , t_n) \mapsto x =
    v + t_1 \alpha_1 + \ldots + t_n \alpha_n ,$$
whose Jacobian is $1/|\Gamma|$, and substituting back into
\eqref{eqAA}, we get
\begin{equation} \labell{cat}
 {\sum_{\bfC \cap \Z^n}}^{\bfq} f = \sum_{\gamma \in \Gamma}
\prod_{i=1}^n {\bf N}_{q_i}^{k,\lambda_{\gamma,i}}(\deldel{h_i})
\left. \int\limits_{\bfC(h)} f(x) dx \right|_{h=0} +
R_{\bfq,k}^{\bfC} (f) ,
\end{equation}
where $\bfC(h_1,\ldots,h_n)$ denotes the image of the ``dilated''
standard orthant ${\bf O}(h_1,\ldots,h_n)$ under the affine
transformation $L$, and where the remainder is given by
\begin{equation}
R^{\bfC}_{\bfq,k}(f) := \frac{1}{|\Gamma|} \sum_{\gamma \in
\Gamma} R^{st}_{\bfq,k}
(\lambda_{\gamma,1},\dots,\lambda_{\gamma,n} ; L^* f),
\end{equation}
where
$$ \lambda_{\gamma,j} := \e^{2 \pi i \left< \gamma, \alpha_j \right> }.$$

\bigskip

Let now $\Delta$ be a simple integral polytope, given by \eqref{Delta}. Choose a
polarizing vector for $\Delta$ and let $\bfC_v^\sharp$ denote the polarized tangent
cones. We can consider the ``dilated polytope" $\Delta(h_1,\ldots,h_d)$, which is
obtained by shifting the $i$th facet outward by a ``distance" $h_i$. More precisely,
$$ \Delta(h) = \bigcap_{i=1}^d \{ x \mid \langle u_i , x \rangle
   + \mu_i + h_i \geq 0 \}
\quad \text{where} \quad h = (h_1,\ldots,h_d). $$

\noindent Then $\Delta(h)$ is simple if $h$ is sufficiently small.
The polar decomposition of $\Delta(h)$ involves ``dilated
orthants". However, dilating the facets of $\Delta$ outward
results in dilating some facets of $\bfC_v^\sharp$ inward and some
outward. Explicitly, for $i \in I_v = \{ i_1, \ldots, i_n \} $,
the inward normal vector to the $i$th facet of $\bfC_v^\sharp$ is
\begin{equation} \labell{u i v sharp}
   u_{i,v}^\sharp = \begin{cases}
 u_i & \text{ if } \alpha_{i,v}^\sharp =  \alpha_{i,v} \\
-u_i & \text{ if } \alpha_{i,v}^\sharp =  -\alpha_{i,v} .
\end{cases}
\end{equation}
Hence, the dilated orthants that occur on the right hand side of
the polar decomposition of $\Delta(h)$ are
$\bfC_v^\sharp(h_{i_1,v}^\sharp, \ldots, h_{i_n,v}^\sharp)$, where
\begin{equation} \labell{h i v sharp}
 h_{i,v}^\sharp = \begin{cases}
 h_i & \text{ if } \alpha_{i,v}^\sharp =  \alpha_{i,v}  \\
-h_i & \text{ if } \alpha_{i,v}^\sharp = - \alpha_{i,v}.
\end{cases}
\end{equation}

Because the inward normals to the facets of $\bfC_v^\sharp$ are
given by \eqref{u i v sharp}, the dual basis to these vectors is
$\alpha_{i,v}^\sharp$, $i \in I_v$, and the roots of unity that
appear in the Euler-Maclaurin formula for $\bfC_v^\sharp$ are then
\begin{equation} \labell{lambda sharp}
\lambda_{\gamma,i,v}^\sharp = \e^{2 \pi i \left< \gamma,
\alpha_{i,v}^\sharp \right> } = \begin{cases}
 \lambda_{\gamma,i,v} & \text{ if } \alpha_{i,v}^\sharp =  \alpha_{i,v} \\
\lambda_{\gamma,i,v}\inv & \text{ if } \alpha_{i,v}^\sharp =
-\alpha_{i,v} .
\end{cases}
\end{equation}

Let $k > 1$ be an integer. For any compactly supported function
$f$ on $\R^n$ of type $C^{nk}$, we then get from \eqref{wheee} and
\eqref{cat}
\begin{multline}\label{prev to key}
{\sum_{\Delta \cap \Z^n}}^{q} f =
\sum_v (-1)^{\# v} {\sum_{\bfC_v^\sharp \cap \Z^n}}^{\bfq^\sharp_v} f \\
 = \sum_v (-1)^{\# v} \sum_{\gamma \in \Gamma_v}
   \prod_{j \in I_v = \{ i_1, \ldots, i_n \} }
   {\bf N}_{q^\sharp_{j,v}}^{k,\lambda_{\gamma,j,v}^\sharp}(\deldel{h_{j,v}^\sharp}) \left.
   \int\limits_{\bfC_v^\sharp(h_{i_1,v}^\sharp,\ldots,h_{i_n,v}^\sharp)}
   f(x) dx \right|_{h=0} \\
   + R_{q,k}^{\Delta}(f),
\end{multline}
where the remainder is given by
\begin{equation} \labell{remdef}
R^\Delta_{q,k} (f) := \sum_v
(-1)^{\#v}R^{C_v^\sharp}_{\bfq^\sharp_v,k} (f)
\end{equation}
and where the $h_{i,v}^\sharp$'s are given in \eqref{h i v sharp}.
Note that either $h_{i,v}^\sharp = h_{i}$,
$\lambda_{\gamma,i,v}^\sharp = \lambda_{\gamma,i,v}$ and
$q_{i,v}^\sharp=q$, or $h_{i,v}^\sharp = -h_{i}$,
$\lambda_{\gamma,i,v}^\sharp = \lambda_{\gamma,i,v}\inv$ and
$q_{i,v}^\sharp=1-q$. By the symmetry property
\eqref{Dandlambdainverse}, we have
\begin{equation}\label{key}
{\bf N}_{q^\sharp_{i,v}}^{k,\lambda_{\gamma,i,v}^\sharp}
(\deldel{h_{i,v}^\sharp})
 = {\bf N}_{q}^{k,\lambda_{\gamma,i,v}} (\deldel{h_i}).
\end{equation}

Now, for $j \not\in I_v$, because $\lambda_{\gamma,j,v} = 1$ (see
~(\ref{lambda is one})), we have $\bfN^{k,\lambda_{\gamma,j,v}
}(\deldel{h_j}) = 1 + $powers of $\deldel{h_j}$. Also for $j \not
\in I_v$, the cone $C_v^\sharp( h_{i_1,v}^\sharp, \ldots,
h_{i_n,v}^\sharp)$ is independent of $h_j$. These facts together
with Formula~(\ref{key}) imply that~(\ref{prev to key}) is equal
to
\begin{equation} \labell{eqA}
\left. \sum_v (-1)^{\# v} \sum_{\gamma \in \Gamma_v}
   \prod_{j=1}^d \bfN^{k,\lambda_{\gamma,j,v}} (\deldel{h_j})
   \int\limits_{C_v^\sharp(h_{i_1,v}^\sharp, \ldots, h_{i_n,v}^\sharp)}
   f(x) dx \right|_{h=0} + R^\Delta_{q,k}(f).
\end{equation}

Because $\lambda_{\gamma,j,F} = \lambda_{\gamma,j,v}$ whenever $v
\in F$ (see~(\ref{horse})), we can define
\begin{equation}
\labell{def of Mk gamma F} {\bf N}^k_{q,\gamma,F} = \prod_{j=1}^d
{\bf N}_q^{k,\lambda_{\gamma,j,F}}(\deldel{h_j}) \quad \text{ for
} \gamma \in \Gamma_F ,
\end{equation}
and we have
\begin{equation} \labell{Fv}
{\bf N}^k_{q,\gamma,F} = {\bf N}^k_{q,\gamma,v} \quad \text{
whenever } \gamma \in \Gamma_F \text{ and } v \in F,
\end{equation}
where we identify $\gamma \in \Gamma_F$ with its image under the
inclusion map $\Gamma_F \hookrightarrow \Gamma_v$.

Then \eqref{eqA} is equal to
\begin{eqnarray}
\nonumber \left.\sum_v (-1)^{\# v} \sum_{\gamma \in \Gamma_v} {\bf
N}_{q,\gamma,v}^k
\int_{C_v^\sharp(h_{i_1}^\sharp,\ldots,h_{i_n}^\sharp)}
    f(x) dx \right|_{h=0} + R^\Delta_{q,k}(f) \\
\labell{eqB}
  =  \left.\sum_F \sum_{\gamma \in \Gamma_F^\flat}
       {\bf N}_{q,\gamma,F}^k \sum_{v \in F} (-1)^{\# v}
       \int_{C_v^\sharp(h_{i_1}^\sharp,\ldots,h_{i_n}^\sharp)}
       f(x) dx \right|_{h=0} + R^\Delta_{q,k}(f) ,
\end{eqnarray}
by~(\ref{Gammas}) and \eqref{Fv}. In the interior summation we may
now add similar summands that correspond to $v \not \in F$. These
summands make a zero contribution to (\ref{eqB}) for the following
reason.  If $v \not \in F$ then there exists $i \in I_F \ssminus
I_v$. Because $i \not \in I_v$, the cone
$C_v^\sharp(h_{i_1,v}^\sharp , \ldots h_{i_n,v}^\sharp)$ is
independent of $h_i$. So it is enough to show that ${\bf
N}_{q,\gamma,F}^k$ is a multiple of $\deldel{h_i}$. But because
$\gamma \in \Gamma_F^\flat$ and $i \in I_F$, we have
$\lambda_{\gamma,i,F} \neq 1$. (See Claim \ref{claim3}.) By Remark
\ref{M vs L}, this implies that ${\bf
N}_q^{k,\lambda_{\gamma,i,F}}(\deldel{h_i})$, which is one of the
factors in ${\bf N}_{q,\gamma,F}^k$, is a multiple of
$\deldel{h_i}$. Hence, \eqref{eqB} is equal to
\begin{multline}
\left.\sum_F \sum_{\gamma \in \Gamma_F^\flat}
       {\bf N}_{q,\gamma,F}^k \sum_{\text{all } v} (-1)^{\# v}
       \int_{C_v^\sharp( h_{i_1,v}^\sharp ,\ldots, h_{i_n,v}^\sharp )}
       f(x) dx \right|_{h=0} + R^\Delta_{q,k}(f)  \\
     = \left.\sum_F \sum_{\gamma \in \Gamma_F^\flat}
       {\bf N}_{q,\gamma,F}^k \int_{\Delta(h)} f(x) dx \right|_{h=0}
       + R^\Delta_{q,k}(f) .\\
\end{multline}

We have therefore proved our main result:

\begin{Theorem}\labell{main}
Let $\Delta$ be a simple integral polytope in $\R^n$. Let $k>1$
and let $f\in C^{nk}_c(\R^n)$ be a compactly supported function on
$\R^n.$  Choose a polarizing vector for $\Delta$.
Then
$${\sum_{\Delta \cap \Z^n}}^{q} f = \left.\sum_F
\sum_{\gamma \in \Gamma_F^\flat}
       {\bf N}_{q,\gamma,F}^k \int_{\Delta(h)} f(x) dx \right|_{h=0}
       + R^\Delta_{q,k}(f) $$
where ${\bf N}_{q,\gamma,F}^k$ are differential operators defined
in \eqref{def of Mk gamma F} and where the remainder
$R^\Delta_{q,k}(f)$ is given by equation (\ref{remdef}). Moreover,
the differential operators ${\bf N}_{q,\gamma,F}^k$ are of order
$\leq k$ in each of the variables $h_1,\ldots,h_d$. Also, the
remainder can be expressed as a sum of integrals over orthants of
bounded periodic functions times various partial derivatives of
$f$ of order no less than $k$ and no more than $kn$. This
remainder is independent of the choice of polarization and is a
distribution supported on the polytope $\Delta$.
\end{Theorem}

In particular, if $q=1/2$, we get Theorem 2 of \cite{ArX}, while
if $q=1$, we have a formula for the ordinary, unweighted sum.
Notice that if $q=0$, we also have an unweighted sum but only over
the interior lattice points in the polytope.

%$${\sum_{\Delta \cap \Z^n}} f = \left.\sum_F
%\sum_{\gamma \in \Gamma_F^\flat}
%       \prod_{i}\todd_{\gamma}^{2k} (\deldel h_i) \int_{\Delta(h)} f(x) dx \right|_{h=0}
%       + R^\Delta_{1,k}(f)$$

% =========================================================
\section{Estimates on the remainder and
         an exact Euler-Maclaurin formula for polynomials}
% =========================================================

In order to derive a formula for polynomials from Theorem
\ref{main}, we first require an estimate on the remainder term
$R^\Delta_{q,k}(f)$. Returning to the definition of functions
$Q_{m,\lambda}(x)$ in \cite{ArX}, we see that
$Q_{m,\lambda}(x)$ is a periodic function on $\R$.  It follows
that $Q_{m,\lambda}(x)$ is bounded.  Since the operators ${\bf
N}_q^{k,\lambda}(\deldel{h})$ are differential operators of order
$k$, the definition of the remainder
$R^{st}_{\bfq,k}(\lambda_1,\dots,\lambda_n;f)$ (see
\eqref{remainder}) shows that

\begin{equation}\labell{orthantest}
|R^{st}_{\bfq,k}(\lambda_1,\dots,\lambda_n;f)| \leq K_{\bfq} \cdot
{\rm sup}_{\{j_1,\dots,j_n\}}
|\partial_1^{j_1}\dots\partial_n^{j_n} f|_{L_1(\R^n)},
\end{equation}

\noindent where the supremum is taken over all $n$-tuples
$\{j_1,\cdots,j_n\}$ with $k \leq j_1 +\cdots + j_n  \leq nk.$

The definition of $R^\Delta_{q,k}$ and equation (\ref{orthantest})
then give the same estimate for the remainder on the polytope, of
course with a different constant.

\begin{Proposition}\labell{polyest}
The remainder term in Theorem \ref{main} can be estimated by

$$|R^\Delta_{q,k}(f)|
\leq K_q(k,\Delta) \cdot {\rm sup}_{\{j_1,\dots,j_n\}}
|\partial_1^{j_1}\dots\partial_n^{j_n} f|_{L_1(\R^n)},$$

\noindent where the supremum is taken over all $n$-tuples
$\{j_1,\cdots,j_n\}$ with $k \leq j_1 + \cdots + j_n  \leq nk.$
\end{Proposition}

The estimate in Proposition \ref{polyest} implies

\begin{Proposition}\labell{polys}
Let $p$ be a polynomial on $\R^n,$ and choose $k \geq {\rm deg}~p
+ n + 1.$ Then
$${\sum_{\Delta \cap \Z^n}}^q p =
\left.\sum_F \sum_{\gamma \in \Gamma_F^\flat}
       {\bf N}_{q,\gamma,F}^k \int_{\Delta(h)} p(x) dx \right|_{h=0}.$$

\end{Proposition}

\begin{Corollary}\label{regular}
Let $p$ be a polynomial and suppose that the polytope $\Delta$ is regular.
Then
$${\sum_{\Delta\cap\Z^n}}^q p = \left.\prod_{i=1}^d {\bf\chi}_q (\deldel h_i)
\int_{\Delta(h)} p(x) dx \right|_{h=0}.$$
\end{Corollary}

\end{document}